\documentclass{article}
\usepackage{amsmath,amsthm}
\usepackage{amssymb}
\usepackage{enumerate}
\usepackage{xypic}

\newtheorem{thm}{Theorem}[section]

\theoremstyle{definition}

\numberwithin{equation}{section}

\newcommand{\bc}{\begin{center}}
\newcommand{\ec}{\end{center}}
\newcommand{\ra}{\rightarrow}
\newcommand{\lra}{\longrightarrow}

\begin{document}

\title{A Natural Occurrence of Shift Equivalence}

\author{Fatma Muazzez \c{S}im\c{s}ir\\
\c{C}\i nar Sokak No:115/2\\
06170 Ankara, Turkey\\
E-mail: fsimsir@gmail.com
\and 
Cem Tezer\\
Department of Mathematics,\\ Middle East Technical University\\
06535 Ankara, Turkey\\
E-mail: rauf@metu.edu.tr}

\date{}

\maketitle


\renewcommand{\thefootnote}{}



\renewcommand{\thefootnote}{\arabic{footnote}}
\setcounter{footnote}{0}


\begin{abstract}
A natural occurrence of  shift equivalence in a purely
algebraic setting   constitutes the subject matter of the
following short exposition.
\end{abstract}

\section{Introduction}
Group endomorphisms  $\alpha : G \longrightarrow G$, $\beta : H \longrightarrow H$,
are said to be {\em ~conjugate~}
if there exists an isomorphism $\theta :G \longrightarrow H$  
such that  $\theta \circ \alpha = \beta \circ \theta .$
$\alpha : G \longrightarrow G$ and $\beta : H \longrightarrow H$ are said
to be {\em ~shift~~equivalent~}  if there exist group endomorphisms
$\varphi : G \longrightarrow H$ , $\psi : H \longrightarrow G$ and $n \in \mathbb Z^+$  such that
the relations\\
\begin{eqnarray*}
   \varphi \circ \alpha   &=& \beta \circ \varphi \\
    \psi \circ \beta      &=& \alpha \circ \psi  \\
    \psi \circ \varphi    &=& \alpha^{n}  \\
     \varphi \circ \psi   &=& \beta^{n}  
\end{eqnarray*}
hold, equivalently, the diagrams\\\\
\xymatrix{& G \ar[d]_\varphi \ar[r]^\alpha & G \ar[d]^\varphi \\
& H \ar[r]^\beta & H } ~~~~ \xymatrix{& G \ar[r]^\alpha & G  \\
& H \ar[u]^\psi  \ar[r]^\beta & H \ar[u]_\psi}  ~~~~ \xymatrix{& G \ar[d]_\varphi \ar[r]^{\alpha^n} & G \ar[d]^\varphi \\
& H \ar[ur]^\psi \ar[r]^{\beta^n} & H }\\\\
commute.This state of affairs is described by 
saying that $\varphi$, $\psi$  effect a
shift equivalence of $\alpha$ to $\beta$ of  lag  ~$n \in \mathbb Z^+$.\\ 
The concept of shift equivalence was introduced 
by R. F. Williams \cite{Williams1}, \cite{Williams2} in the context of topological
dynamics. The fact that shift equivalence is an equivalence relation
among group endomorphisms  can be demostrated by a straightforward
argument. \cite{Tezer1}\\
Clearly both conjugacy and 
shift equivalence can be defined in any 
category and the former constitutes a 
special case of the latter in two ways: 
\begin{itemize}
\item  A shift equivalence with lag 0
is a conjugacy.\\ 
\item  A shift equivalence between
two automorphisms is a conjugacy. 
\end{itemize}
The simple  result presented here was independently observed by
Yu. I. Ustinov in a short report \cite{Ustinov}. In our opinion
this is the most straightforward and  natural occurrence of shift equivalence as complete invariant.
Although by no means entirely novel,
we feel that this elegant result 
deserves to be  available to a wider public
in the form of an independent exposition.\\
Another very natural occurrence 
of shift equivalence 
arises in shape and homotopy theory, \cite{Tezer2}.
\section{Statement and proof of the main result:}
Given a group endomorphism  $\alpha : G \lra G$ the
{\em ~simple direct limit~}  of $\alpha$   ,
denoted by
$\mathfrak G$ = $ \lim_{\ra}(G, \alpha)$, is the set of equivalence
classes in $ G \times \mathbb Z^+$  under the equivalence 
relation $\sim$  where
\[   (g,n)\sim (g',n') \]
 iff
\[  \alpha^{N - n}(g)  =\alpha^{N - n'}(g')  \]
for some $N \geq n, n'$.

$\sim$  can be easily checked to be an equivalence 
relation. $\mathfrak G$  has a natural group
structure with respect to the binary operation
\[
           (g,n)(g',n') = (\alpha^{n'} (g) \alpha^{n} (g'), n+n')
\]
where,  by abuse of notation,  we 
let $(g,n)$   stand for the equivalence 
class it represents.
Again  it can be routinely 
checked that this is a well-defined operator satisfying all
group axioms. There are  two    natural isomorphisms on $\mathfrak G$ : Firstly,
\[
\check{\alpha} : \mathfrak G \lra \mathfrak G
\]
defined by
\[
\check{\alpha} ((g,n)) = (\alpha (g), n)
\]
secondly,
\[
s_{\alpha} : \mathfrak G  \lra \mathfrak G
\]
(which we like to call the  ``coshift'')  defined by
\[
s_{\alpha} ((g,n)) = (g,n + 1).
\]
Again, well-definedness and morphology need 
checking. To see that  $\check{\alpha}$  and
$s_{\alpha}$ are isomorphisms it is enough to observe 
that
\[ 
\check{\alpha} \circ s_{\alpha} = 
s_{\alpha} \circ \check{\alpha} = Id(\mathfrak G).
\]

\begin{thm}
Let $G$ and $H$ be 
finitely generated groups ,
$\alpha : G \lra G$,
   $\beta  : H \lra H$  group endomorphisms, 
$\mathfrak G$ = $\lim_{\ra} (G, \alpha )$ , $\mathfrak H$ = 
$\lim_{\ra} (H , \beta )$. The isomorphisms $s_{\alpha}$ : 
$\mathfrak G$ $\lra$
$\mathfrak G$,  $s_{\beta}$ : $\mathfrak H$ $\lra$ $\mathfrak H$  are conjugate iff $\alpha$,  $\beta$  are
shift equivalent.
\end{thm}

\begin{proof}  Given a subset $K$ of  
a group,  
let $\langle K \rangle$  denote the
subgroup generated by $K$.
There exist finite sets $A \subseteq G$,  
$B \subseteq H$   such that $G =\langle A\rangle$
 ,\   $H = \langle B \rangle$.
Assume first that $s_{\alpha}$ and  
$s_{\beta}$ , or equivalently,  $\check{\alpha}$ and
$\check{\beta}$ are conjugate: There exists an isomorphism
\[
T : \mathfrak G  \lra \mathfrak H
\]
such that
\[
T \circ \check{\alpha}  = \check{\beta} \circ T.
\]
Let
\[
i_{\alpha} : G \lra \mathfrak G
\]
\[
i_{\beta} : H \lra  \mathfrak H
\]
be the natural injections defined by
\[
i_{\alpha} (g) = (g,0) \in \mathfrak G
\]
\[
i_{\beta} (h) = (h,0) \in \mathfrak H.
\]
We have
\[
T \circ i_{\alpha} (G) \subseteq \langle T \circ i_{\alpha} (A) \rangle.
\]
Clearly $T \circ i_{\alpha} (A)$  is a finite subset of $\mathfrak H$.  Hence there exists $k \in \bf Z^{+}$
such that
\[
    T \circ i_{\alpha} (G) \subseteq \langle T \circ i_{\alpha}(A) \rangle \subseteq H \times \{k\}.
\]
Therefore,
\[
   \check{\beta}^{k} \circ T \circ i_{\alpha}(G) \subseteq H \times \{0\}.
\]
We define
\[
\varphi = i_{\beta}^{-1} \circ 
\check{\beta}^{k} \circ T \circ i_{\alpha} : G \lra H.
\]
Similarly, there exists a  sufficiently large $l \in \bf Z_{+}$  such that
\[
\psi =i_{\alpha}^{-1} \circ 
\check{\alpha}^{l} \circ T \circ i_{\alpha} : H \lra G
\]
is a well defined homomorphism. We claim  that $\varphi$ and 
$\psi$ effect a
shift equivalence of $\alpha$  to  $\beta$  with 
lag $ k + l \in \bf Z^{+}$ : Clearly
\[
     \varphi \circ \alpha = \beta \circ \varphi
\]
\[
    \psi \circ \beta =\alpha \circ \psi.
\]
Moreover,
\[
\psi \circ \varphi = i_{\alpha}^{-1} \circ \check{\alpha}^{l} \circ T^{-1} \circ i_{\beta} \circ
i_{\beta}^{-1} \circ \check{\beta}^{k} \circ T \circ i_{\alpha} = \alpha^{k+l}.
\]
Similarly,
\[
\varphi \circ \psi = \beta^{k+l}.
\]
Conversely assume,  that there exist $\varphi :G \lra H$ ,
$\psi : H \lra G$  and $n \in \bf Z^{+}$
such that   $\varphi \circ \alpha = \beta \circ \phi$  ,
$\psi \circ \beta = \alpha \circ \psi$  ,
$\psi \circ \varphi = \alpha^{n}$ ,
 $\phi \circ \psi  = \beta^{n}$. Consider the map
\[
E : \mathfrak G \lra \mathfrak H
\]
defined by
\[
    E((g,m)) = (\varphi(g),m)
\]
and note that $E$ is well-defined: If  $\alpha^{l-m}(g) 
= \alpha^{l-m'}(g')$ ,  then
\[
\varphi \circ \alpha^{l-m}(g) = \varphi \circ \alpha^{l-m'}(g')
\]
hence
\[
\beta^{l-m} \circ \varphi (g) = \beta^{l-m'} \circ \varphi(g').
\]
We have also
        $E \circ \check{\alpha} = \check{\beta} \circ E$
owing to  $\varphi \circ \alpha  = \beta \circ \varphi$ ,  once
again. Similarly define
\[
F : \mathfrak H \lra \mathfrak G
\]
by
\[
      F((h,m)) = (\psi(h),  m).
\]
We observe
\[
   F \circ E ((g,m)) = F((\varphi(g),m)) =(\psi \circ \varphi(g),m)
\]
\[
     =(\alpha^{n}(g),m)
\]
\[
     =\check{\alpha}^{n}(g,m).
\]
Thus $F \circ E = \check{\alpha}^{n}$ .  The right 
hand side is an isomorphism,
$E$ is an isomorphism,  too,  which commutes with 
$\check{\alpha}$  and  $\check{\beta}$ .
\end{proof}

\end{document}